\ifx\control\undefined

\documentclass[10pt, leqno, twoside]{article}
\usepackage{amssymb, amsmath} 
\usepackage{latexsym} 
\usepackage{graphics} 
\usepackage{url} 
\usepackage{enumerate}
\usepackage[all]{xy}
\usepackage{graphicx}
\usepackage{float}
\usepackage{amsrefs}
\usepackage{color}

\usepackage{authblk}
\usepackage{blindtext}

\usepackage{enumitem}

\numberwithin{equation}{section}

\newtheorem{theo}{Theorem}

\newtheorem{prop}[theo]{Proposition}

\newtheorem{prob}[theo]{Problem}

\newenvironment{proo}[1][\proofname]{\normalfont{\itshape
#1{:}}\quad\mdseries\ignorespaces}
{{$\Box$}{\vskip\belowdisplayskip}}
\newcommand{\proofname}{Proof}
\newtheorem{defi}[theo]{Definition}
\newtheorem{rem}[theo]{Remark}
\begin{document}

\fi
%%%%%%%%%%%%%%%%%%%%%%%%%%%%%%%%%%%%%%%%%%%%%%%%%%%%%%

\title{Geometric continuity of plane curves in terms
of Riordan matrices and an application to the
F-chordal problem}

%\author{L. Felipe Prieto-Mart\'inez, Raquel S\'anchez-Cauce}

\author{L. Felipe Prieto-Mart\'inez\thanks{Department of Mathematics, Universidad Aut\'onoma de Madrid (Spain), luisfelipe.prieto@uam.es}}
%\author{}
%\affil[1]{Department of Mathematics, Universidad Aut\'onoma de Madrid (Spain)}
%\affil[2]{Department of Artificial Intelligence, Universidad Nacional de Educaci\'on a Distancia (Spain)}

\date{\today}

\maketitle

\maketitle

\begin{abstract} The first goal of this article is to provide an statement of the conditions for geometric continuity of order $k$, referred in the bibliography as beta-constraints, in terms of Riordan matrices. The second one is to see this new formulation in action to solve a theoretical cuestion about uniqueness of analytic solution for a general and classical problem in plane geometry: the $F$-chordal problem.
\end{abstract}

\section{Introduction}
\label{intro}

 For $k\geq 1$ a plane curve $c$ is of class $G^k$ (has geometric continuity of order $k$) if there exists a local regular parametrization of class $C^k$ of this curve in a neighbourhood of each point of $c$. The case $k=\infty$ can also be considered and then if, in addition, we impose the local regular parametrizations to be analytic, then we say that the curve is analytic.

Geometric continuity is a concept of great importance in computer-aided geometric design, where the objects are frequently described in terms of parametric splines. For more information see, for instance, the articles \cite{B.GC1, B.GC2} or the book \cite{K.GC}. 

One of the main problems related to geometric continuity of plane curves can be stated as follows:

\begin{prob}\label{problem.GC} Let $c$ be a curve with a (countinuous) piecewise defined parametrization $\gamma:(-\varepsilon,\varepsilon)\to c$ given by
$$\gamma(t)=\begin{cases}\gamma_{left}(t)=(x_{left}(t),y_{left}(t)) & t\in(-\varepsilon,0)\\
V & t=0\\ \gamma_{right}(t)=(x_{right}(t),y_{right}(t)) & t\in (0,\varepsilon) \end{cases} $$

\noindent where $\gamma_{left},\gamma_{right}$ are parametrizations of class $C^k$.  We will call the point of intersection $V=\gamma(0)$ the vertex.

Despite the fact of $\gamma_{left},\gamma_{right}$ being $C^k$, a regular parametrization of $c$ may not exist in any neighbourhood of $V$ (see Figure \ref{fig.leftright}). Assume that the Taylor polynomials of degree $k$ of the functions $x_{left}(t)$, $y_{left}(t)$, $x_{right}(t)$, $y_{right}(t)$ at $t=0$ exist. Which compatibility conditions should satisfy the coefficients of these Taylor polynomials if we want $c$ to be a curve of class $G^k$?

\begin{figure}
\centering 
  \includegraphics[width=70mm]{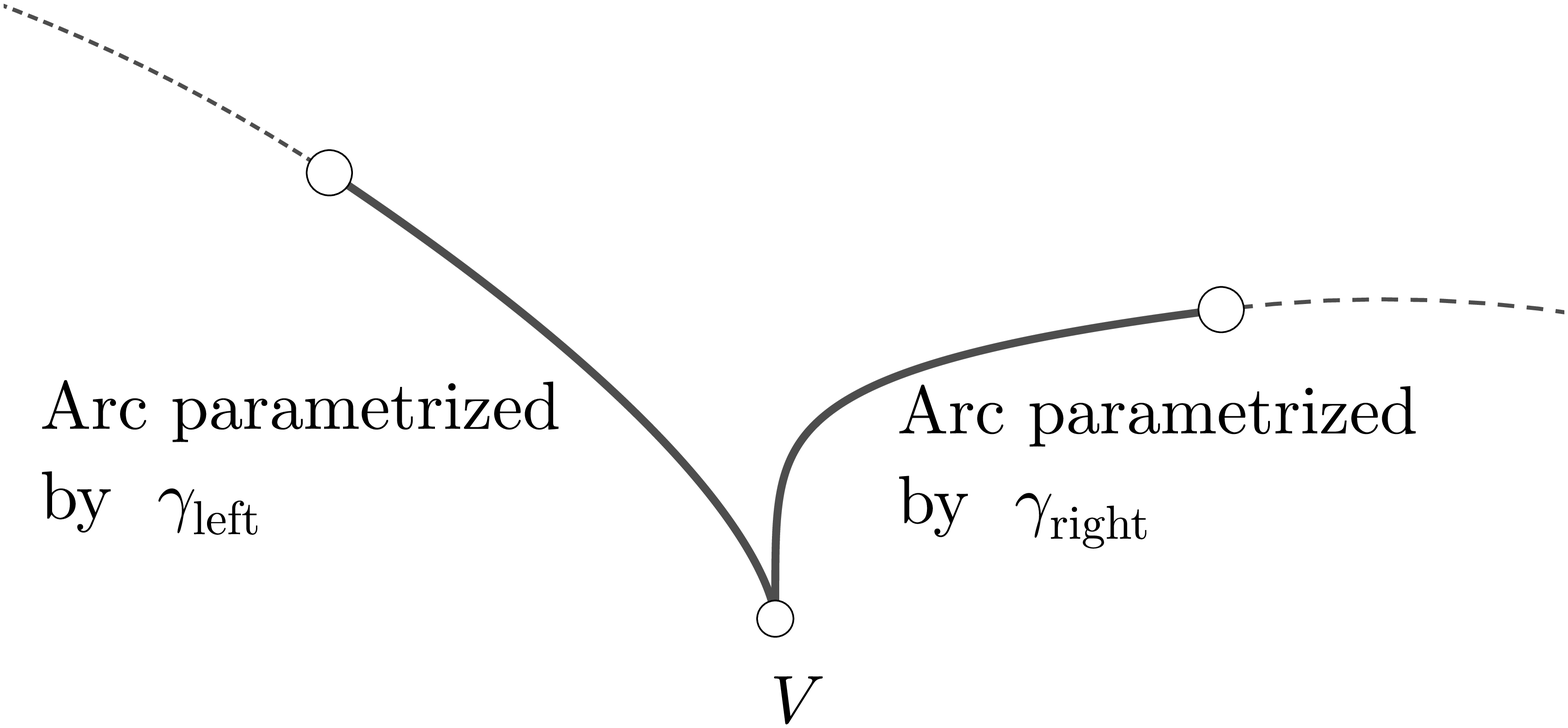}
\label{fig.leftright}       % Give a unique label
\end{figure}

\end{prob}

As explained briefly in the abstract, this article has two main goals. The first one is to provide a new and  useful formulation of the conditions on the Taylor polynomial required in Problem \ref{problem.GC}. This target is reached in Section \ref{sect.GCRM}. The conditions for the Taylor polynomial  are sometimes expressed in the bibliography in terms of the so called \emph{connection matrices} (see section 2.1 in \cite{K.GC}). The  statement presented here (Remark \ref{rem.statGC}) is done in terms of \emph{Riordan matrices}. Finite and infinite Riordan matrices have a well studied structure and properties (some basics are provided in Section \ref{sect.riordan}) which are more adequate for doing a certain kind of computations. For example, the inverse limit structure studied in \cite{LMMPS.IL}, allow us to do easily proofs by induction,  like the ones required in our second target. This second goal is to show in action this new formulation solving a theoretical problem. We show how this new statement can be used to solve a questions about uniqueness of analytic solutions appearing in a well known problem in plane geometry (Problem \ref{problem.interior}).

Let $F$ be a symmetric function in two variables, defined in $[0,\infty)\times [0,\infty)$. Given a convex region $B$ which boundary is $c$, a chord in $c$ is any segment which endpoints belong to $c$. We say that $P$ in the interior of $B$ is an \emph{interior $F$-chordal point} if there exists a constant $K_P$ such that for every chord in $c$, passing through $P$ and with endpoints $A,B$ we have $F(\left|A-P\right|,\left|B-P\right|)=K_P$.

 \begin{figure}[!h]
  \centering
    \includegraphics[width=120mm]{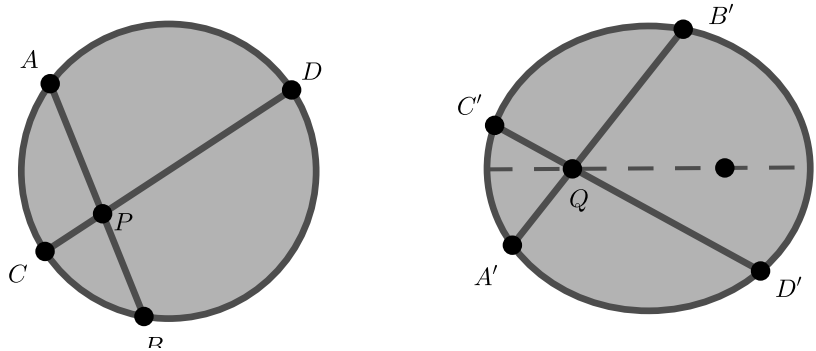}
\caption{The center of a disk is an F-chordal point for $F(a,b)=a+b$. Any point in the disk is an $F$-chordal point for $F(a,b)=a\cdot b$ as a consequence of \emph{Steiner's Power of a Point Theorem} (so $\left|P-A\right|\cdot \left|P-B\right|=\left|P-C\right|\cdot \left|P-D\right|$). Klee in \cite{K} noticed that any of the two focus of an ellipse is an F-chordal point for $F=\frac{1}{a}+\frac{1}{b}$ (so $\frac{1}{\left|P-A'\right|}+\frac{1}{ \left|P-B'\right|}=\frac{1}{\left|P-C'\right|}+\frac{1}{ \left|P-D'\right|}$)}
  \label{fig.fpoint}
\end{figure}

\begin{prob}[Two points Interior F-chordal Problem] \label{problem.interior} For a given symmetric function in two variables $F$ defined in $[0,\infty)\times[0,\infty)$, find a plane Jordan curve which interior region is convex, with two interior $F$-chordal points.

\end{prob}

For any solution $c$ of this problem, we call the line through $P,Q$ the \emph{axis}, and the two points where the axis meet $c$ the \emph{vertices}. As we will see later, one of this vertices will also be a \emph{vertex} in the sense of Problem \ref{problem.GC}.

The $F$-chordal Problem is a generalization proposed in the book \cite{CFG} of an older problem: the  \emph{Equichordal Problem}, stated in 1916-1917 independently by Fujiwara \cite{Fuji} and Blaschke, Rothe and Weitzenb\"ock \cite{BRW}. For a convex region $B$ with boundary $c$, a point $P$ in the interior of $B$ is an \emph{equichordal point} if all the  chords of $c$ passing through $P$ are of the same length (the center of a circle is an equichordal point, for instance). The Equichordal Problem ask wether a convex region $B$  can have  two equichordal points. In \cite{Fuji}, the author already proved that no such a region can have three or more equichordal points. But we had to wait until 1997 when Rychlik \cite{R} answered in a negative way the question. In the meantime, Wirsing showed in \cite{W} that, if such a region exists then the curve $c$ of the boundary should be analytic. This is one of the reasons why analytic solutions for Problem \ref{problem.interior} are of interest. Other problems (some of which remain open)  with an interesting history can be considered as  particular cases of the $F$-chordal Problem too. More is said about this in Section 6.

In Section \ref{sect.G1}, $G^1$ solutions of Problem \ref{problem.interior} are considered. In Section \ref{sect.GT} we prove the following theorem, which answer the question of uniqueness of analytic solutions of Problem \ref{problem.interior} for most of the cases studied in the bibliography.

\begin{theo} \label{th.main} Let four different collinear points $V_{1},P,Q,V_{2}$, where $P,Q$ are between $V_{1},V_{2}$. Let $F:[0,\infty)^2\to \mathbb R$ be a  symmetric function in two variables, satisfying the following conditions:

\begin{enumerate}

\item[(i)]  it is $C^\infty$,

\item[(ii)]  $\frac{\partial F}{\partial b}\mid_{(\|P-V_1\|,\|P-V_2\|)}, \frac{\partial F}{\partial b}\mid_{(\|Q-V_1\|,\|Q-V_2\|)}\neq 0$,

\item[(iii)]  $\frac{\partial F}{\partial a}\mid_{(\|P-V_1\|,\|P-V_2\|)},\frac{\partial F}{\partial a}\mid_{(\|Q-V_1\|,\|Q-V_2\|)} \neq 0 $. 

\item[(iv)] $\forall n\in\mathbb N$, 
$$\left(\frac{\left.\frac{\partial F}{\partial a}\right|_{(\|Q-V_1\|,\|Q-V_2\|)}\cdot \left.\frac{\partial F}{\partial b}\right|_{(\|P-V_1\|,\|P-V_2\|)}}{\left.\frac{\partial F}{\partial a}\right|_{(\|P-V_1\|,\|P-V_2\|)}\cdot \left.\frac{\partial F}{\partial b}\right|_{(\|Q-V_1\|,\|Q-V_2\|)}}\right)^n\neq  \frac{\|V_{2}-Q\|}{\|V_{2}-P\|}\frac{\|V_{1}-P\|}{\|V_{1}-Q\|} $$

\end{enumerate}

\noindent Then if there exists an analytic solution for the interior $F$-chordal Problem with  interior $F$-chordal points $P,Q$ and vertices $V_1,V_2$, this solution is unique.

\end{theo}

 In this case, $k_P=F(\|V_2-P\|,\|V_1-P\|)$ and $k_Q=F(\|V_2-Q\|,\|V_1-Q\|)$.

As a consequence of this theorem, we will provide an alternative proof of part of the results involved in the articles by Rychklik \cite{R, R2} for the Equichordal Problem. We also comment on some of the most studied cases appearing in the bibliography, and we give a generalization of Theorem \ref{th.main} for a different definition of $F$-chordal point that does  not require the curve $c$ to be a Jordan cuve. This is done in Section 6.

Finally, in Section \ref{section.future}, we propose some problems of uniqueness of analytic solutions in plane geometry, that are suitable to be studied using the techniques in this article .

\section{Basics on Riordan matrices} \label{sect.riordan}

Riordan matrices first appeared in \cite{S.O}, although the original definition was slightly different to the one in current use. The classical survey \cite{S.S} contains more information with a similar notation to the one used here. Riordan matrices and generalized Riordan matrices are special types of infinite lower triangular matrices:
$$(a_{ij})_{i,j=0}^\infty=\begin{bmatrix} a_{00} \\ a_{10} & a_{11} \\ a_{20} & a_{21} & a_{22} \\ \vdots & \vdots & \vdots & \ddots\end{bmatrix} $$

\begin{defi} An infinite matrix $(a_{ij})_{i,j=0}^\infty$ is a \emph{generalized Riordan matrix} over the reals if and only if  there exists two formal  power series:
\begin{equation} \label{eq.defriordan} d\in\mathbb R[[t]], \qquad\qquad\qquad h\in t\mathbb R[[t]] \end{equation}

\noindent such that, for every $0\leq i,j$, $a_{ij}=[t^i](d\cdot h^i)$, where $[t^i]f$ denotes the $i-th$ coefficient of the formal power series $f$. In other words, the generating function of the  $i$-th column $ a_{0i}, a_{1i}, a_{2i}, \hdots $ is $d\cdot h^i$. In this case, we write: $(a_{ij})_{i,j=0}^\infty=R(d,h)$. If we replace the condition \eqref{eq.defriordan} by the stronger one $d\in\mathbb R[[t]]\setminus t\mathbb R[[t]]$, $ h\in t\mathbb R[[t]]\setminus t^2\mathbb R[[t]]$, we would have an \emph{(ordinary) Riordan matrix} instead.

\end{defi}

\noindent The condition $h\in t\mathbb R[[t]]$ ensures that  gereneralized Riordan matrices are always lower triangular.

\medskip

Generalized Riordan matrices have an important property known as the \emph{First Fundamental Theorem of Riordan Matrices} (1FTRM). If we multiply any generalized Riordan matrix by an infinite column vector, we obtain a new infinite column vector:
$$R(d,h)\begin{bmatrix} F_0\\ F_1 \\ F_2 \\ \vdots \end{bmatrix} =\begin{bmatrix} G_0\\ G_1 \\ G_2 \\ \vdots \end{bmatrix} $$

\noindent And if $F$ is the generating function of the sequence $F_0,F_1,F_2,\ldots $, then the generating function of $G_0, G_1, G_2,\ldots$ is $d\cdot (F\circ h)$.  As a consequence of this theorem, for every two generalized Riordan matrices, we have that
\begin{equation}R(d,h)\cdot R(f,g)=R(d\cdot f\circ h,g\circ h) \end{equation}

\noindent It can be proved straightforward that the set of (ordinary) Riordan matrices is a group, which a description of the inverse in terms of the corresponding formal power series too. But this is not necessary for this article.

\medskip

We will need something else, concerning the \emph{inverse limit structure} of the Riordan group. This structure will allow us to do proofs by induction. Define a generalized partial Riordan matrix $R_n(d,h)$ to be the principal submatrix of size $(n+1)\times (n+1)$ of a generalized Riordan matrix $R(d,h)=(a_{ij})_{0\leq i,j <\infty}$
$$R(d,h)=\left[\begin{array}{c c c | c c} & & & 0 & \hdots \\ & R_n(d,h)& & \vdots &  \\ & & & 0 & \hdots  \\ \hline a_{n+1,0} & \hdots & a_{n+1,n} & a_{n+1,n+1} \\ \vdots & & \vdots & \vdots &\ddots  \end{array}\right] $$

\begin{rem} \label{rem.riordan} 
\noindent As a consequence of matrix block multiplication for triangular matrices, we have that:
$$R(d,h)\begin{bmatrix} F_0 \\ F_1 \\ F_2 \\ \vdots \end{bmatrix} =\begin{bmatrix} G_0 \\ G_1 \\ G_2 \\ \vdots \end{bmatrix}\qquad \Longleftrightarrow \qquad \forall  n\in\mathbb N,\quad R_n(d,h)\begin{bmatrix} F_0 \\ F_1 \\ \vdots \\ F_n \end{bmatrix} =\begin{bmatrix} G_0 \\ G_1 \\ \vdots \\ G_n \end{bmatrix} $$

Moreover, see that in the matrix $R_n(d,h)$ depends only on the coefficients of  $Taylor_n(d),Taylor_n(h)$ ($Taylor_n(f)$ denotes the Taylor polynomial of degree $n$ at $t=0$ of $f$). In particular,
$$R_n(d,h)=R_n(\widetilde d,\widetilde h)\Longleftrightarrow \begin{cases} Taylor_n(d)=Taylor_n(\widetilde d)\\ Taylor_n(h)=Taylor_n(\widetilde h) \end{cases}$$

See also that if
\begin{equation}\label{eq.riordaninduction} R_n(d,h)\begin{bmatrix} F_0 \\ F_1 \\ \vdots \\ F_n \end{bmatrix} =\begin{bmatrix} G_0 \\ G_1 \\ \vdots \\ G_n \end{bmatrix}\end{equation}

\noindent holds for $n=k$, then it holds for every $m\leq k$.

Aditionally, if for $n=k$ we already have a partial Riordan matrix $R_k(d,h)$ satisfying \eqref{eq.riordaninduction} then we can search for a matrix $R_{k+1}(d,h)$ satisfying \eqref{eq.riordaninduction} for $n=k+1$ just looking at the last entry in the column vector obtained in each side of the equality, and we have only two new parameters in $R_{k+1}(d,h)$ to achieve this. We will call this proccess \emph{extending} the matrix.

\end{rem}

Much more can be said about generalized partial Riordan matrices. For example an intrinsic definition (not depending on the definition of a ``bigger'' generalized Riordan matrix) is also possible. We recommend  \cite{LMMPS.IL} for more information about this finite dimensional matrices.

\section{Geometric Continuity in terms of Riordan matrices} \label{sect.GCRM}

Now that we have introduced Riordan matrices, we will go back to Problem \ref{problem.GC}. In the context of this problem, $c$ is $G^k$ if there exists a $C^k$ regular reparametrization $u$ of $\gamma_{left}$ (or equivalently of $\gamma_{right}$) such that $u(0)=0$ and 
$$\widetilde\gamma:(-\delta,\delta)\longrightarrow  c $$
$$\widetilde\gamma(t)=\begin{cases} \gamma_{left}(u(t)) & t\in(-\delta,0)\\ V &  t=0\\ \gamma_{right}(t) & t\in(0,\delta)\end{cases} $$

\noindent  For each $k$, let  the corresponding Taylor polynomials at $t=0$
$$Taylor_n(x_{left}(t))=a_0+a_1t+\ldots+a_kt^k,\;Taylor_n(y_{left}(t))=b_0+b_1t+\ldots+b_kt^k $$
$$Taylor_n(x_{right}(t))=c_0+c_1t+\ldots+c_kt^k,\;Taylor_n(y_{right}(t))=d_0+d_1t+\ldots+d_kt^k $$

\noindent  The conditions that this parameters $a_i,b_i,c_i,d_i$ for $i=0,\ldots,k$ must satisfy for the curve $c$ are a set of linear equations. In the bibliography, the linear relations (equivalent to those proposed here) are known  as \emph{beta-constraints}, and is frequently stated in terms of the so called \emph{connection matrices} (see \cite{K.GC}). But we suggest here to express these conditions in terms of Riordan matrices, which provide us a powerful tool for doing computations.

\begin{rem} \label{rem.statGC} The conditions for having geometric continuity of order $n$ at $V$ in the notation above are equivalent to the existence of a partial (ordinary) Riordan matrix $R_k(1,u)$ such that:
\begin{equation} \label{eq.new} R_n(1,u)\begin{bmatrix}a_0\\ \vdots \\ a_n \end{bmatrix}=\begin{bmatrix}c_0\\ \vdots \\ c_n \end{bmatrix}\qquad\qquad R_n(1,u)\begin{bmatrix}b_0\\ \vdots \\ b_n \end{bmatrix}=\begin{bmatrix}d_0\\ \vdots \\ d_n \end{bmatrix}\end{equation}

\noindent Recall that $u\in t\mathbb R[[t]]\setminus t^2\mathbb R[[t]]$.

\end{rem}

This statement is equivalent to the classic one involving the \emph{connection matrices}. But it has two advantages: (1) If, for instance, $(x_{left},y_{left})$ are fixed and known functions, we hace a bridge between functional equations and the linear restrictions for the corresponding Taylor polynomials (the same occurs for different types of restrictions between  $(x_{left},y_{left})$ and  $(x_{right},y_{right})$). (2) For the $G^\infty$ curves case, we can easily do proofs by induction (as we will do to show uniqueness of analytic solutions for Problem \ref{problem.interior}).

In this second case, we need to find, $\forall n\in\mathbb N$, matrices $R_n(1,u)$ satisfying \eqref{eq.new}. The stategy (which is the one used  in Theorem \ref{th.main}) is the following. First we find a solution for the case $n=2$. Then we \emph{extend} (as explained in Remark \ref{rem.riordan}) this solution to a solution for the case $n=3$ and so on. if we have a matrix $R_k(1,u)$ satisfying \eqref{eq.new} for $n=k$, we then prove that there is a unique matrix $R_{k+1}(1,u)$ satisfying this same equation for $n=k+1$. The existence of this matrices often implies certain restrictions between the coefficients $a_i,b_i,c_i,d_i$.

\section{$G^1$ solutions for Problem \ref{problem.interior}} \label{sect.G1}

First of all, we will discuss the necessity of the conditions imposed for the function $F$ in Theorem \ref{th.main}. Usually, for most of the particular cases appearing in the bibliography of the $F$-chordal Problem, in the equation $F(a,b)=k_P$  we can write $b$ explicitely as a function $\varphi_P(a)$. For example:

\begin{itemize}

\item For the \emph{Equichordal Problem}, $F(a,b)=a+b=k_P$, so we can consider that $b=k_P-a$ .

\item For the \emph{Equiproduct Problem} (except in the trivial case $k_P=0$), $F(a,b)=a\cdot b=k_P$, and again we can take $b=\frac{k_P}{a}$ .

\end{itemize}

\noindent But, what can we do in the rest of the cases?
 
\begin{rem} For an arbitrary function in two variables $F(a,b)$, the Implicit Function Theorem, together with condition (ii)  implies that in a neighbourhood of $(\|P-V_1\|,\|P-V_2\|),(\|Q-V_1\|,\|Q-V_2\|)$ we can find two real functions $\varphi_P,\varphi_Q$ such that
$$F(|A-P\|,\|B-P\|)=k_P\Longleftrightarrow \|B-P\|=\varphi_P(\|A-P\|) $$
$$F(|A-Q\|,\|B-Q\|)=k_Q\Longleftrightarrow  \|B-Q\|=\varphi_Q(\|A-Q\|)$$

\noindent Moreover, condition (i) ensures that $\varphi_P,\varphi_Q$ are $C^\infty$ functions. We will denote the corresponding Taylor series as:
$$\varphi_P(a)=\varphi_{P0}+\varphi_{P1}(a-\|P-V_{1}\|)+\varphi_{P2}(a-\|P-V_{1}\|)^2+\ldots$$
$$\varphi_Q(a)=\varphi_{Q0}+\varphi_{Q1}(a-\|Q-V_{1}\|)+\varphi_{Q2}(a-\|P-V_{1}\|)^2+\ldots$$

 Now we can re-state conditions $(iii),(iv)$ in terms of these functions $\varphi_P,\varphi_Q$. See that:
$$\varphi_{P1}=-\frac{\frac{\partial F}{\partial a}\mid_{(\|P-V_{1}\|,\|P-V_{2}\|)}}{\frac{\partial F}{\partial b}\mid_{(\|P-V_{1}\|,\|P-V_{2}\|)}},\varphi_{Q1}=-\frac{\frac{\partial F}{\partial a}\mid_{(\|Q-V_{1}\|,\|Q-V_{2}\|)}}{\frac{\partial F}{\partial b}\mid_{(\|Q-V_{1}\|,\|Q-V_{2}\|)}}, $$

\noindent so condition $(iii)$ implies $\varphi_{P1},\varphi_{Q1}\neq 0$, and condition $(iv)$ implies that 
\begin{equation} \label{eq.(vi)}\forall n\in\mathbb N,\qquad \left(\frac{\varphi_{Q1}}{\varphi_{P1}}\right)^n\neq \frac{\|V_{2}-Q\|}{\|V_{2}-P\|}\frac{\|V_{1}-P\|}{\|V_{1}-Q\|} =\frac{\varphi_{Q0}}{\varphi_{P0}}\frac{x_0-1}{x_0+1} \end{equation}

\end{rem}

The relation between Problem \ref{problem.interior} and Problem \ref{problem.GC} arises from the fact that each $F$-chordal point  induces an involutive correspondence between points in $c$. The image of a point $A$ through this correspondence is the other point lying in the intersection between $c$ and the line through the $F$-chordal point and $A$. And thank to this, we have that:

\begin{rem}  \label{rem.rep} Let $c$ be a curve with two $F$-chordal points $P,Q$. Any parametrization $\gamma:(-\varepsilon,\varepsilon) \to c$ with $\gamma(t)=(x(t),y(t))$ of $c$ in a neighbourhood of one of the vertices $V_1$ induces two (one for each $F$-chordal point) parametrizations of $c$ in a neighbourhood of $V_2$.
$$\gamma_P(t)= P+\varphi_P(|P-(x(t),y(t))\|)\frac{P-(x(t),y(t))}{\|P-(x(t),y(t))\|}$$
$$\gamma_Q(t)= Q+\varphi_Q(|Q-(x(t),y(t))\|)\frac{Q-(x(t),y(t))}{\|Q-(x(t),y(t))\|}$$

 \begin{figure}[!h]
  \centering
    \includegraphics[width=70mm]{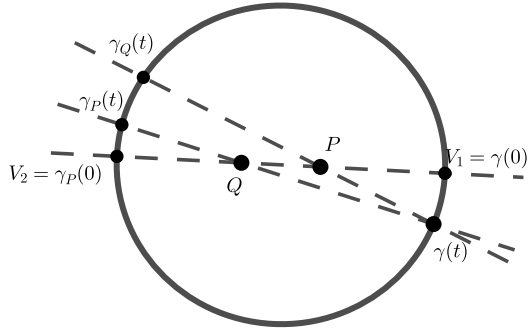}
\caption{This picture may help to understand this double parametrization and the correspondence between points in $c$.}
  %\label{F.I}
\end{figure}

So if we want a solution of Problem \ref{problem.interior} to be $G^n$ near $V_2$,  it should exist a $C^n$ regular reparametrization $u(t)$, such that $u(0)=0$ and such that satisfies
\begin{equation}\label{eq.rpr1} \gamma_P(u(t))=\gamma_Q(t)\end{equation}

\end{rem}

Equation \eqref{eq.rpr1} is a restrictions for the coefficients of the Taylor polynomial of degree $n$ of $x(t),y(t)$. And we can state \eqref{eq.rpr1} as a functional equation, suitable to be expressed in terms of Riordan matrices, as done in Remark \ref{rem.statGC}. This leads to the following:

\begin{rem} From now on, to study \eqref{eq.rpr1}, we will take $P=(1,0)$, $Q=(-1,0)$, $V_{1}=(x_0,0)$ with $x_0>1$.

In the notation of the previous remark, to find a $G^n$ solution for Problem \ref{problem.interior}, we need to solve $\forall k\leq n$ the following system of two matricial equations:

\begin{equation} \label{eq.mainn}\begin{scriptsize} \begin{cases}\begin{bmatrix}  2 \\ 0 \\ \vdots \\ 0 \end{bmatrix} +R_n(1,u)R_n(1-x,\sqrt{(1-x)^2+y^2}-(x_0-1))\begin{bmatrix} \uparrow\\  Coefficients \; of\\ \frac{\varphi_P((x_0-1)+t)}{x_0-1+t} \\ \downarrow \end{bmatrix} =\\
\qquad\qquad\qquad= R_n(-1-x,\sqrt{(x+1)^2+y^2}-(x_0+1))\begin{bmatrix} \uparrow\\  Coefficients \; of\\ \frac{\varphi_Q((x_0+1)+t)}{x_0+1+t} \\ \downarrow \end{bmatrix}\\
R_n(1,u)R_n(y,\sqrt{(1-x)^2+y^2}-(x_0-1))\begin{bmatrix} \uparrow\\  Coefficients \; of\\ \frac{\varphi_P((x_0-1)+t)}{x_0-1+t} \\ \downarrow \end{bmatrix} =\\
\qquad\qquad\qquad= R_n(y,\sqrt{(x+1)^2+y^2}-(x_0+1))\begin{bmatrix} \uparrow\\  Coefficients \; of\\ \frac{\varphi_Q((x_0+1)+t)}{x_0+1+t} \\ \downarrow \end{bmatrix} \end{cases} \end{scriptsize}\end{equation}

\end{rem}

Then we have the following:

\begin{prop} \label{lemm.base} Let $c$ be $G^1$ solution for Problem \ref{problem.interior}, with vertices $V_{1},V_{2}$, interior $F$-chordal points $P,Q$, and such that $F$ satisfies conditions (i)-(iv) from Theorem \ref{th.main}. Then the tangent vector of $c$ at $V_{1}$ is either parallel, either perpendicular to the axis of $c$. 
\end{prop}

\begin{proo} Suppose that we have a parametrization $\gamma(t)=(x(t),y(t))$ such that $\gamma(0)=V_{1}=(x_0,0)$ and such that the corresponding Taylor polynomials at $t=0$ are
$$Taylor_1(x)(t)=x_0+x_1t \qquad\qquad Taylor_1(y)(t)=y_1t $$

\noindent  Then the case $n=1$ of \eqref{eq.mainn} is
$$\begin{cases}\begin{bmatrix}  2 \\ 0 \end{bmatrix} +\begin{bmatrix}  1 \\ 0 & u_1 \end{bmatrix}\begin{bmatrix} 1-x_0 \\ -x_1 & (1-x_0)x_1\end{bmatrix}\begin{bmatrix}\frac{\varphi_{P0}}{x_0-1} \\ -\frac{\varphi_{P0}}{(x_0-1)^2}+\frac{\varphi_{P1}}{(x_0-1)}\end{bmatrix} =\\
\qquad\qquad\qquad\qquad= \begin{bmatrix} -1-x_0 \\ -x_1 & (-1-x_0)x_1\end{bmatrix}\begin{bmatrix}\frac{\varphi_{Q0}}{x_0+1} \\ -\frac{\varphi_{Q0}}{(x_0+1)^2}+\frac{\varphi_{Q1}}{(x_0+1)}\end{bmatrix}\\
\;\\
\begin{bmatrix}  1 \\ 0 & u_1 \end{bmatrix}\begin{bmatrix} 0 \\ y_1 & 0\end{bmatrix}\begin{bmatrix}\frac{\varphi_{P0}}{x_0-1} \\ -\frac{\varphi_{P0}}{(x_0-1)^2}+\frac{\varphi_{P1}}{(x_0-1)}\end{bmatrix}  = \begin{bmatrix} 0 \\ y_1 & 0\end{bmatrix}\begin{bmatrix}\frac{\varphi_{Q0}}{x_0+1} \\ -\frac{\varphi_{Q0}}{(x_0+1)^2}+\frac{\varphi_{Q1}}{(x_0+1)}\end{bmatrix}\end{cases}$$

\noindent where $Taylor_n(u)(t)=u_1t$. This system leads to a system of 4 equations, each of them corresponding to one of the entries of the column vectors of length two obtained in each side of each matricial equation:
$$\begin{cases} 2-\varphi_{P0}=-\varphi_{Q0}\\ 
-u_1\varphi_{P1}x_1=-\varphi_{Q1}x_1\\
0=0\\
u_1y_1\frac{\varphi_{P0}}{x_0-1}=y_1\frac{\varphi_{Q0}}{x_0+1}\end{cases}$$

\noindent The first and third equations are trivial, so we only need to matter the other two. Taking into account that $\gamma(t)$ is a regular parametrization and so $\gamma'(t) =(x_1,y_1)\neq (0,0)$, we have three possible cases:

\begin{itemize}

\item Case 1: $x_1,y_1\neq 0$, which is not possible since implies a contradiction with condition $(iv)$ (see Equation \eqref{eq.(vi)}, which recall that is a version of condition (iv) in terms of $\varphi_P,\varphi_Q$)
$$\begin{cases}u_1\varphi_{P1}=\varphi_{Q1}\\
u_1\frac{\varphi_{P0}}{x_0-1}=\frac{\varphi_{Q0}}{x_0+1}\end{cases}\Rightarrow \frac{\varphi_{Q1}}{\varphi_{P1}}=\frac{\varphi_{Q0}}{\varphi_{P0}}\frac{x_0-1}{x_0+1}$$

\item Case 2: $x_1\neq 0,y_1=0$
$$\begin{cases}u_1\varphi_{P1}=\varphi_{Q1}\\
0=0\end{cases}\Rightarrow u_1=\frac{\varphi_{Q1}}{\varphi_{P1}}$$

\item Case 3: $x_1=0, y_1\neq 0$
$$\begin{cases}0=0\\
u_1\frac{\varphi_{P0}}{x_0-1}=\frac{\varphi_{Q0}}{x_0+1}\end{cases}\Rightarrow u_1=\frac{\varphi_{Q0}}{\varphi_{P0}}\frac{x_0-1}{x_0+1}$$

\end{itemize}

\hspace{10cm}\end{proo}

\section{Proof of Theorem \ref{th.main}}\label{sect.GT}

First of all, we need to check that Case 2 in the end of the previous proof does not correspond to any analytic solution of Problem \eqref{problem.interior}.

\begin{prop} \label{prop.line} Let $F$ satisfying conditions (i)-(iv) of Theorem \ref{th.main}. An analytic  curve $c$ that satisfies \eqref{eq.mainn} for every $n\in\mathbb N$ is a line segment contained in the line through $P,Q$. So it does not correspond to any  analytic solution of Problem \ref{problem.interior}, which should be a Jordan curve with  $P,Q$ in its interior region.

\end{prop}

\begin{proo} We are going to do the proof by induction. This argument is similar but more simple that the one of Theorem \ref{th.main}. The base case has already been considered in Proposition \ref{lemm.base}. Let
$$x_0+x_1t+\ldots+x_{k+1}t^{k+1}\qquad\qquad y_2t^2+\ldots+y_{k+1}t^{k+1} $$

\noindent be the Taylor polynomials at $t=0$ of $x(t),y(t)$ respectively. Assume that, for some $k\geq 1$, if $y_0,\ldots,y_k=0$ and we have some $x_0,\ldots,x_k$, $u_1,\ldots,u_k$ that are a solution of \eqref{eq.mainn} for $n=k$. Then the second equation of the system \eqref{eq.mainn} for $n=k+1$  is of the type:
$$\begin{bmatrix}1 \\ 0 & u_1 \\ \vdots &  \vdots & \ddots \\ 0& u_{k+1} & \hdots & u_{1}^{k+1} \end{bmatrix} \begin{bmatrix}0  \\ \vdots  & \ddots \\ 0 &\hdots& 0 \\ y_{k+1} & 0 & \hdots & 0   \end{bmatrix} \begin{bmatrix} \uparrow\\  Coefficients \; of\\ \frac{\varphi_P((x_0-1)+t)}{x_0-1+t} \\ \downarrow \end{bmatrix}=$$
$$\qquad\qquad\qquad= \begin{bmatrix}0  \\ \vdots  & \ddots \\ 0 &\hdots& 0 \\ y_{k+1} & 0 & \hdots & 0   \end{bmatrix}\begin{bmatrix} \uparrow\\  Coefficients \; of\\ \frac{\varphi_Q((x_0+1)+t)}{x_0+1+t} \\ \downarrow \end{bmatrix}$$

\noindent and we are going to see that it implies that $y_{k+1}=0$.

 The result in each side of the matricial equation is a column vector. The equality between the two last entries in each is:
$$u_{1}^{n+1}y_{n+1} \frac{\varphi_{P0}}{x_0-1}=y_{n+1}\frac{\varphi_{Q0}}{x_0+1}\Rightarrow (u_{1}^{n+1}\frac{\varphi_{P0}}{x_0-1}-\frac{\varphi_{Q0}}{x_0+1})y_{n+1}=0$$

\noindent The only solution of the above equation is $y_{n+1}=0$, since $u_1=\frac{\varphi_{Q1}}{\varphi_{P1}}$ and so:
$$u_1^{n+1}\neq \frac{\varphi_{Q0}}{\varphi_{P0}}\frac{x_0-1}{x_0+1} $$

\noindent This shows that the case $x_1\neq 0, y_1=0$ implies the Taylor series of $y(t)$ at 0 equals 0. And so, in a neighbourhood of $V_1$, $c$ is a segment contained in the line through $P,Q$. Since $c$ is analytic, by the \emph{Principle of Analytic Continuation} $c$ must be a line segment, which cannot be  a solution for Problem \ref{problem.interior}.

\hspace{10cm}\end{proo}

Now that we have discard Case 2, we now that any analytic solution for Problem \ref{problem.interior} has its tangent vector at any of its vertices perpendicular to its axis, and we can complete the proof of Theorem \ref{th.main}.

\medskip

\begin{proo}[Proof of Theorem \ref{th.main}] Let the Taylor series of $x(t),y(t)$ (recall that $\gamma(t)=(x(t),y(t))$ is a parametrization of an analytic solution $c$ in a neighbourhood of $V_1$) be respectively
$$x(t)=x_0+x_2t^2+\ldots\qquad\qquad y(t)=y_1t+y_2t^2+\ldots $$

 We know (Proposition \ref{lemm.base}) that, since the tangent vector of $c$ at $V_1$ is perpendicular to the axis, then \eqref{eq.mainn} for $n=1$ implies that  $\displaystyle{u_1=\frac{\varphi_{Q0}}{\varphi_{P0}}\frac{x_0-1}{x_0+1}}$.

We are going to prove by induction in $k$ the following statement: for each choice of $x_0>1$ (this value is determined by the vertex $V_1$), $y_1\neq 0$, and  $y_2,\ldots,y_n$, there exists a unique choice of $x_2,\ldots,x_k$, $u_1,\ldots,u_k$ such that \eqref{eq.mainn} holds. This determines univocally $\gamma(t)$ up to reparametrization (this is the reason of the freedom of the parameters $y_2,\ldots,y_n,\ldots)$ and thus $c$, according to the \emph{Principle of Analytic Continuation}.

Although we have already studied the case $k=1$, for a better understanding of this proof, we include the case $k=2$ which is the first significant one. In this case \eqref{eq.mainn} is:
$$\begin{cases}\begin{bmatrix}  2 \\ 0 \\ 0 \end{bmatrix} +\begin{bmatrix}  1 \\ 0 & u_1 \\ 0 & u_2 & u_1^2 \end{bmatrix}\begin{bmatrix} 1-x_0 \\ 0 & 0\\-x_2& \frac{1}{2}[ (x_0-1)x_2-y_1^2] & 0\end{bmatrix}\begin{bmatrix}\frac{\varphi_{P0}}{x_0-1} \\ a_1\\ a_2 \end{bmatrix} = \\
\qquad\qquad\qquad =\begin{bmatrix} -1-x_0 \\ 0 & 0  \\ -x_2 & \frac{|}{2}[(x_0+1)x_2+y_1^2] & 0 \end{bmatrix}\begin{bmatrix}\frac{\varphi_{Q0}}{x_0+1} \\b_1 \\ b_2\end{bmatrix}\\
\;\\
\begin{bmatrix}  1 \\ 0 & u_1 \\ 0 & u_2 & u_1^2 \end{bmatrix}\begin{bmatrix} 0 \\ y_1 & 0 \\ y_2 & 0 & 0 \end{bmatrix}\begin{bmatrix}\frac{\varphi_{P0}}{x_0-1} \\a_1 \\ a_2 \end{bmatrix}  = \begin{bmatrix} 0 \\ y_1 & 0 \\ y_2 & 0 & 0 \end{bmatrix}\begin{bmatrix}\frac{\varphi_{Q0}}{x_0+1} \\b_1 \\ b_2 \end{bmatrix}\end{cases}$$

\noindent where $\frac{\varphi_P((x_0-1)+t)}{x_0-1+t}=\frac{\varphi_{P0}}{x_0-1}+a_1x+a_2x^2+\ldots$ and $\frac{\varphi_P((x_0+1)+t)}{x_0+1+t}=\frac{\varphi_{Q0}}{x_0+1}+b_1x+b_2x^2+\ldots$

\medskip

\noindent Each matricial equation lead to a system of 3 linear equations in the indeterminates $x_1,x_2,u_1,u_2$:
\begin{equation} \label{eq.proof3} \begin{cases} 2-\varphi_{P0}=-\varphi_{Q0}\\ 
-u_1\varphi_{P1}x_1=-\varphi_{Q1}x_1\\
-u_1^2\frac{\varphi_{P0}}{x_0-1}x_2+[\frac{1}{2}u_1^2a_1((x_0-1)x_2-y_1^2)]=\\
\qquad=-\frac{\varphi_{Q0}}{x_0+1}x_2+[\frac{1}{2}b_1((x_0+1)x_2+y_1^2)]\end{cases} \end{equation}
\begin{equation} \label{eq.proof32} \begin{cases} 0=0\\
u_1y_1\frac{\varphi_{P0}}{x_0-1}=y_1\frac{\varphi_{Q0}}{x_0+1}\\ 
\frac{\varphi_{P0}}{x_0-1}y_1u_2+[u_1^2\frac{\varphi_{P0}}{x_0-1}y_2]=[\frac{\varphi_{Q0}}{x_0+1}y_2]\end{cases}\end{equation}

\noindent We have already discussed in Proposition \ref{lemm.base} the values of $u_1$ that make the two first equations in each system to hold. The last equation in \eqref{eq.proof32} does not depend on $x_2$, and has nontrivial coefficient of the indeterminate $u_2$ (the coefficient is $\displaystyle{y_1\frac{\varphi_{P0}}{x_0-1}}$), so it has a unique solution in this indeterminate. On the other hand, the last equation in \eqref{eq.proof3} has again a non-trivial coefficient (according to the hypothesis of the theorem and the value of $u_1$) for the indeterminate $x_2$ (the coefficient is $\displaystyle{-\frac{\varphi_{P0}}{x_0-1}u_1^2+\frac{\varphi_{Q0}}{x_0+1}}$) and so it has a unique solution in the indeterminate $x_2$ too.

Now assume that the statement  is true for $k-1\geq 2$. We want to solve \eqref{eq.mainn}, for $n=k$. The second matricial equation in \eqref{eq.mainn} is of the type:
$$ \begin{bmatrix}1 \\ 0 & u_1 \\ \vdots & \vdots & \ddots \\ 0 & u_n & \hdots & u_1^k \end{bmatrix}\begin{bmatrix}0 \\ y_1 & 0 \\ \vdots & \ddots & \ddots \\  y_k & \hdots & \ddots & 0 \end{bmatrix}\begin{bmatrix}\frac{\varphi_{P0}}{x_0-1}\\ a_1 \\ \vdots \\ a_k \end{bmatrix}=\begin{bmatrix}0 \\ y_1 & 0 \\ \vdots & \ddots & \ddots \\  y_k & \hdots & \ddots & 0 \end{bmatrix}\begin{bmatrix}\frac{\varphi_{Q0}}{x_0+1}\\ b_1 \\ \vdots \\ b_k  \end{bmatrix}$$

\noindent The equation corresponding to the last entry in the column vector of each side, is a linear equation in the indeterminate $u_k$ of the form: 
$$y_1\frac{\varphi_{P0}}{x_0-1}u_n+[C_1]=[C_2]$$

\noindent where nothing in the brackets depend on $u_k, x_k$ (they do on $y_k$). So we have a unique solution on the indeterminate $u_k$ that makes this equation hold. On the other hand, the first matricial equation in \eqref{eq.mainn} is of the type:
$$\begin{bmatrix} 2 \\ 0 \\ \vdots \\ 0 \end{bmatrix} + \begin{bmatrix}1 \\ 0 & u_1 \\ \vdots & \vdots & \ddots \\ 0 & u_k & \hdots & u_1^n \end{bmatrix}\begin{bmatrix}1-x_0 \\ 0 & 0 \\ \vdots & \ddots & \ddots \\  -x_k & \hdots & \ddots & 0 \end{bmatrix}\begin{bmatrix}\frac{\varphi_{P0}}{x_0-1}\\ a_1 \\ \vdots \\ a_n  \end{bmatrix}=\begin{bmatrix}-1-x_0 \\ 0 & 0 \\ \vdots & \ddots & \ddots \\  -x_k & \hdots & \ddots & 0 \end{bmatrix}\begin{bmatrix}\frac{\varphi_{Q0}}{x_0+1}\\ b_1 \\ \vdots \\ b_k  \end{bmatrix}$$

\noindent The equation corresponding to the last entry in the column vector is:
$$-u_1^k\frac{\varphi_{P0}}{x_0-1}x_k+[C_3]=-\frac{\varphi_{Q0}}{x_0+1}x_k+[C_4] $$

\noindent And nothing in the brackets depends on $x_k$ (they do on $y_k,u_k$). The number ${\frac{\varphi_{Q0}}{x_0+1}-u_1k\frac{\varphi_{P0}}{x_0-1}} $ is not zero according to the hypothesis of the theorem. So this equation has a unique solution in the indeterminate $x_k$.

\hspace{10cm}\end{proo}

\section{Consequences of Theorem \ref{th.main}} \label{sect.con}

First of all, we want to point out that in the proof of Theorem \ref{th.main} we have not used the fact that $c$ must be the boundary of a convex region, neither a Jordan curve. This theorem still holds for a more general definition of interior $F$-chordal point:

\begin{defi} We say that $P$ is an \emph{interior $F$-chordal point} of a curve $c$ if there exists $k_P$ such that for every line through $P$ either (1) $l$ does not intersect $c$ or (2) $l$ meets $c$ at two points $A,B$ satisfying $F(\|A-P\|,\|B-P\|)=k_P$ and such that $P$ is in the interior of the segment $AB$. 

\end{defi}

\noindent For example, for $F=a-b$, the center of symmetry $P$ of any hyperbola is  an interior $F$-chordal point in this sense, with $k_P=0$. With this definition, if we have two point $P,Q$ the line through them can also be considered to be a solution of Problem \ref{problem.interior} (modifying the statement of Proposition \ref{prop.line}).

Secondly, we want to remark that, with almost the same proof, we can obtain this more general version of Theorem \ref{th.main} that will be required in this section:

\begin{theo}\label{th.mainint*}  Let four different collinear points $V_{1},P,Q,V_{2}$, where $P,Q$ are between $V_{1},V_{2}$. Let $F:[0,\infty)^2\to \mathbb R$ be a  symmetric function in two variables, satisfying the following conditions:

\begin{enumerate}

\item[(i)]  it is $C^\infty$,

\item[(ii)]  $\frac{\partial F}{\partial b}\mid_{(\|P-V_1\|,\|P-V_2\|)}, \frac{\partial F}{\partial b}\mid_{(\|Q-V_1\|,\|Q-V_2\|)}\neq 0$,

\item[(iii)]  $\frac{\partial F}{\partial a}\mid_{(\|P-V_1\|,\|P-V_2\|)},\frac{\partial F}{\partial a}\mid_{(\|Q-V_1\|,\|Q-V_2\|)} \neq 0 $. 

\end{enumerate}

\begin{enumerate}

\item[(iv*)] $\forall n\in\mathbb N$, $n\geq 2$, 
$$\left(\frac{\left.\frac{\partial F}{\partial a}\right|_{(\|Q-V_1\|,\|Q-V_2\|)}\cdot \left.\frac{\partial F}{\partial b}\right|_{(\|P-V_1\|,\|P-V_2\|)}}{\left.\frac{\partial F}{\partial a}\right|_{(\|P-V_1\|,\|P-V_2\|)}\cdot \left.\frac{\partial F}{\partial b}\right|_{(\|Q-V_1\|,\|Q-V_2\|)}}\right)^n\neq  \frac{\|V_{2}-Q\|}{\|V_{2}-P\|}\frac{\|V_{1}-P\|}{\|V_{1}-Q\|} $$

\end{enumerate}

\noindent Then if we fix the tangent direction at $V_1$, if there exists an analytic solution for the interior $F$-chordal Problem with  interior $F$-chordal points $P,Q$, vertices $V_{1},V_{2}$ and this given tangent direction  at $V_{1}$, this solution must be unique.

\end{theo}

\begin{proo} In Proposition \ref{lemm.base} the case $x_1,y_1\neq 0$ cannot be excluded, but anyway $u_1=\frac{\varphi_{Q0}}{\varphi_{P0}}\frac{x_0-1}{x_0+1}$. Later in the proof of Theorem \ref{th.main}, one of the matrices in the system \eqref{eq.mainn} is different but it does not affect the argument.

\hspace{10cm}\end{proo}

Finally we are going to collect some consequences of theorems \ref{th.main}, \ref{th.mainint*}. Several particular cases of  $F$-chordal points have been studied in the bibliografphy. We have already discussed about  equichordal points ($F(a,b)=a+b$). In Figure \ref{fig.fpoint} we can see on the left a disk, for which any point in the interior is an \emph{equiproduct point} ($F(a,b)=a\cdot b$, see \cite{CFG, F.E, G, Z} for more information) and on the right an ellipse, for which any of its two focus is an \emph{equireciprocal point} ($F(a,b)=\frac{1}{a}+\frac{1}{b}$, see \cite{K}). In general, the family of $F$-chordal points for $F(a,b)=a^\alpha+b^\alpha$ has also be considered, for $\alpha\in\mathbb R$ (see \cite{CFG}).

The following three results are a direct consequence of Theorem \ref{th.main}:

\begin{theo}[concerning the Equichordal Problem] \label{th.equi} For every four collinear points $V_{1}$, $P$, $Q$, $V_{2}$, if it exists an analytic curve with equichordal points $P,Q$ and vertices $V_{1},V_{2}$, this curve is unique.

\end{theo}

\begin{theo}[concerning the Equireciprocal Problem] For every four collinear points $V_{1}$, $P$, $Q$, $V_{2}$, if it exists an analytic curve with two equireciprocal points, this curve is unique. 

If $\|V_{1}-Q\|=\|V_{2}-P\|$, the  ellipse with foci $P,Q$ and major axis the segment with endpoints $V_{1},V_{2}$ is this unique curve (see \cite{F.E,K} to see that such an ellipse has these properties).

\end{theo}

\begin{theo}[concerning the $F$-chordal Problem for $F(a,b)=a^\alpha+b^\alpha$]

For every four collinear points $V_{left},P,Q,V_{right}$, for $F(a,b)=a^\alpha+b^\alpha$, $\alpha\neq 0$, if it exists an analytic curve with two $F$-chordal points, this curve is unique.

\end{theo}

\noindent And this last one is a consequence of Theorem \ref{th.mainint*}:

\begin{theo}[conerning the Equiproduct Problem] \label{th.circles} For every four collinear points $V_{1}$, $P$, $Q$, $V_{2}$, the circles that pass through $V_1,V_2$ are the unique analytic curves with $F$-chordal points $P,Q$ and vertices $V_{1},V_{2}$.

\end{theo}

In relation to Theorem \ref{th.equi}, the fact that it does not contradict the Theorem by Rychlik in \cite{R} deserves a little explanation. The author already discussed in \cite{R2} that the Equichordal Problem had a local analytic regular solution, pointing out the Helfenstein was wrong in his article \cite{H}. From the local point of view, the family of all the interior $F$-chordal problems studied here behaves in a similar way: we have a unique candidate (up to reparametrization) for the power series of a parametrization near a vertex. To study wether those local solutions can be extended or not to solutions of Problem \ref{problem.interior} needs other type of global techniques. For example, for the \emph{Equichordal Problem}, Rychlik proved that extension cannot be done becouse of the hiperbolicity of the problem. But this is not the case of the Equireciprocal Problem, for instance. The techniques used in this article are not suitable for this global analysis. Anyway, maybe the desciption of the coefficients of the local solution (specially the first terms) appearing in the proof of Theorem \ref{th.equi} could be use in the search for a more simple proof of the result by Rychlik, which remains open for the classical statement of the problem (where the solutions must be the boundary of a convex region).

\section{Further work} \label{section.future}

Finally we will state two other classical problems in plane geometry. The techniques appearing in this article seem to be suitable to prove uniqueness of analytic solutions for them, but an improved argument may be required. Those questions are open, up to our knowledge. 
\medskip

There exists another version of the $F$-chordal Problem, that we could call the \emph{Exterior $F$-chordal Problem} which could be stated as follows. For a convex region $B$ with boundary $c$, we say that a point $P$ in the exterior region of $c$ is an exterior $F$-chordal point if there exists a constant $k_P$ such that for every chord with endpoints $A,B$ in $c$, $F(\left|A-P\right|,\left|B-P\right\|)=k_P$.

\begin{prob}[one point and two points exterior $F$-chordal Problems] For a given symmetric function in two variables $F$ defined in $[0,\infty)\times[0,\infty)$:

\begin{enumerate}

\item[(a)]  find a plane Jordan curve which interior region is convex with one exterior F-chordal point,

\item[(b)] find a plane Jordan curve which interior region is convex with two exterior $F$-chordal points.

\end{enumerate}

\end{prob}

\begin{figure}
\centering  
  \includegraphics[width=70mm]{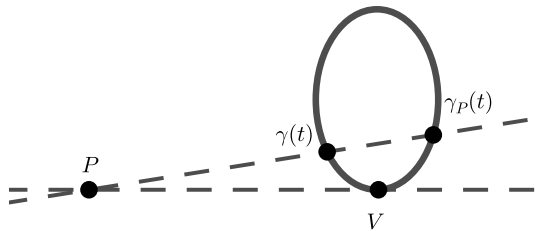}
\caption{Any parametrization $\gamma(t)$ of a neighbourhood of the vertex $V$ (any of the two points where the line from $P$ is tangent to $c$) induces another parametrization $\gamma_P(t)$ in the same neighbourhood}
\label{fig.exterior}       % Give a unique label
\end{figure}

\noindent Some results are known concerning particular cases of this problem, see for example \cite{Z} for the \emph{Equiproduct Problem} case. For these problems,  the double parametrization, analogous to the one for the Interior $F$-chordal Problem that appears in Remark \ref{rem.rep}, can be obtained with only one exterior $F$-chordal point (we omit the  details, but we offer a picture, see Figure \ref{fig.exterior}). So it makes sense to study using our techniques the One Point case. Again, a new definition is possible for \emph{exterior $F$-chordal point}, not requiring $c$ to be a Jordan curve.

On the other hand, we have a problem related to Geometric Tomography, which is a field that focuses on problems of reconstructing plane regions from tomographic data. The term was introduced by R. J. Gardner in the book \cite{G}.  We could state one of the main and most simple problems in this field as follows. 

\begin{prob}\textbf{(One Point, Two Points, One Line or Two Lines Tomographic Reconstruction Problem)} 

A \emph{tomographic image from a point} $P$ in the exterior of a convex region $B$ is a real function $f_P$ such that $f(\theta)$ is the length of   $l_\theta\cap B$, where  $l_\theta$ is the line passing through $P$ and with angle $\theta$ with respect the axis $OX$.   Given one or two tomographic images from a point, find a Jordan curve $c$ such that the interior region of $c$ is convex and has this or these tomographic images.

Equivalently,  the tomographic image from the $OY$ axis $r$,  is a real function $g_r$ such that $g_r(t)$ is the length of the segment  $l_t\cap B$, where $l_t$ the horizontal line which $y$ coordinate equals $t$. The analogue can be defined for any line with the corresponding modifications. Given one or two tomographic images from a line, find a Jordan curve $c$ such that the interior region of $c$ is convex and has this or these tomographic images.

\end{prob}

See \cite{G, GK, V} for more information about these problems. A single tomographic image, either from a point or from a line, induces a double parametrization near a vertex, with a similar picture to the one in Figure \ref{fig.exterior}. Moreover, our method could provide algorithms to approximate the boundary of the region $B$, assuming that it is $G^k$ for some $k$.

%\begin{acknowledgements}
%\end{acknowledgements}




\end{document}